\documentclass[11pt,reqno,a4paper]{amsart}%
\usepackage[a4paper,margin=3cm]{geometry}%
\usepackage[usenames,dvipsnames]{xcolor}%
\usepackage[pagebackref]{hyperref}%
\usepackage{caption,subcaption}%
\usepackage{lmodern,mathtools}%
\usepackage{color,todonotes}%
\usepackage[utf8]{inputenc}%
\usepackage[T1]{fontenc}%
\usepackage{listings}%

\title[Simulating Coulomb and log-gases with HMC]%
{Simulating Coulomb and log-gases\\with Hybrid Monte Carlo algorithms}

\author{Djalil Chafaï}%
\address[DC]{Université Paris-Dauphine, PSL, CNRS, CEREMADE, %
  F-75016 Paris, France}
\email{\url{mailto:djalil(at)chafai.net}}%
\urladdr{\url{http://djalil.chafai.net/}}

\author{Grégoire Ferré}%
\address[GF]{Université Paris-Est, CERMICS (ENPC), INRIA, %
  F-77455 Marne-la-Vallée, France}%
\email{\url{mailto:gregoire.ferre@enpc.fr}}%
\urladdr{\url{https://team.inria.fr/matherials/team-members/gregoire-ferre/}}

\date{Summer 2018, revised Autumn 2018, compiled \today}

\lstdefinelanguage{Julia}%
  {morekeywords={abstract,break,case,catch,const,continue,do,else,elseif,%
      end,export,false,for,function,immutable,import,importall,if,in,%
      macro,module,otherwise,quote,return,switch,true,try,type,typealias,%
      using,while},%
   sensitive=true,%
   morecomment=[l]\#,%
   morecomment=[n]{\#=}{=\#},%
   morestring=[s]{"}{"},%
   morestring=[m]{'}{'},%
}[keywords,comments,strings]%

\newtheorem{theorem}{Theorem}[section]%
\newtheorem{remark}[theorem]{Remark}%

\newtheorem{algo}[theorem]{Algorithm}

\newcommand{\dC}{\mathbb{C}}
\newcommand{\dE}{\mathbb{E}}

\newcommand{\dP}{\mathbb{P}}
\newcommand{\dR}{\mathbb{R}}

\newcommand{\cE}{\mathcal{E}}

\newcommand{\cN}{\mathcal{N}}

\newcommand{\al}{\alpha}
\newcommand{\be}{\beta}
\newcommand{\De}{\Delta}
\newcommand{\de}{\delta}
\newcommand{\ga}{\gamma}

\newcommand{\na}{\nabla}

\newcommand{\Dt}{\Delta t}

\newcommand{\ABS}[1]{{{\left| #1 \right|}}} 

\newcommand{\IND}{\mathbf{1}}

\newcommand{\dd}{\mathrm{d}}
\newcommand{\e}{\mathrm{e}}

\DeclarePairedDelimiter\ceil{\lceil}{\rceil}

\makeatletter
\def\@MRExtract#1 #2!{#1}     
\renewcommand{\MR}[1]{
  \xdef\@MRSTRIP{\@MRExtract#1 !}%
  \href{http://www.ams.org/mathscinet-getitem?mr=\@MRSTRIP}{MR-\@MRSTRIP}}
\makeatother

\numberwithin{equation}{section}



\begin{document}

\begin{abstract}
  Coulomb and log-gases are exchangeable singular Boltzmann--Gibbs measures
  appearing in mathematical physics at many places, in particular in random
  matrix theory. We explore experimentally an efficient numerical method for
  simulating such gases. It is an instance of the Hybrid or Hamiltonian Monte
  Carlo algorithm, in other words a Metropolis--Hastings algorithm with
  proposals produced by a kinetic or underdamped Langevin dynamics. This
  algorithm has excellent numerical behavior despite the singular interaction,
  in particular when the number of particles gets large. It is more efficient
  than the well known overdamped version previously used for such problems,
  and allows new numerical explorations. It suggests for instance to
  conjecture a universality of the Gumbel fluctuation at the edge of beta
  Ginibre ensembles for all beta.
\end{abstract}

\maketitle

{\footnotesize\tableofcontents}

We explore the numerical simulation of Coulomb gases and log-gases by mean of
Hybrid or Hamiltonian Monte Carlo algorithms
(HMC)~\cite{duane-hmc,horowitz1991generalized}. Such algorithms consist
basically in using discretized kinetic (underdamped) Langevin dynamics to
produce proposals for Metropolis--Hastings algorithms. This can be viewed as a
way to add momentum to a Monte Carlo interacting particle system. The basic
outcome of this exploratory work is that HMC algorithms have remarkably good
numerical behavior for such gases despite the singularity of the interactions.
Such algorithms scale well with the dimension of the system,
see~\cite{MR3129023,bou2017geometric}. They are therefore more efficient than
the tamed overdamped version already explored in the literature for instance
in ~\cite{MR3124976}. In this paper, we benchmark the capability of the
algorithm to reproduce known results efficiently, and we make it ready to
explore new conjectures.

Another advantage of this approach is that it could be adapted to take into
account a sub-manifold constraint~\cite{MR2945148}. For
instance, this could be used for simulating random matrices with prescribed
trace or determinant, which is difficult to achieve by direct sampling of
matrices.

For the sake of completeness, we should mention that there are remarkable
alternative simulation algorithms which are not based on a diffusion process,
such as the ones based on piecewise deterministic Markov processes (PDMP), see
for instance \cite{Kapfer-2016} and \cite{PDMPs}.

\section{Boltzmann--Gibbs measures}

We are interested in interacting particle systems subject to an external field
and experiencing singular pair interactions. In order to encompass Coulomb
gases as well as log-gases from random theory, we introduce a vector subspace
$S$ of dimension $d$ of $\mathbb{R}^n$, with $n\geq2$ and $n\geq d\geq 1$. The
particles belong to $S$, and $\mathbb{R}^n$ is understood as a physical
ambient space. We equip $S$ with the trace of the Lebesgue measure of
$\mathbb{R}^n$, denoted by $\dd x$. The external field and the pair
interaction are respectively denoted by $V:S\mapsto\mathbb{R}$ and
$W:S\mapsto(-\infty,+\infty]$, and belong to $\mathcal{C}^2$ functions, with
$W(x)<\infty$ for all $x\neq 0$. For any $N\geq2$, we consider the probability
measure $P_N$ on $S^N=S\times\cdots\times S$ defined by
\begin{equation}
\label{eq:gibbs}
P_N(\dd x)=\frac{\e^{-\be_NH_N(x_1,\ldots,x_N)}}{Z_N}\dd x_1\cdots\dd x_N,
\end{equation}
where $\beta_N>0$ is a parameter,
\[
  Z_N=\int_{S^N}\e^{-\be_NH_N(x_1,\ldots,x_N)}\dd x_1\cdots\dd x_N
\]
is the normalizing factor, and
\[
  H_N(x_1,\ldots,x_N)=\frac{1}{N}\sum_{i=1}^NV(x_i)
  +\frac{1}{2N^2}\sum_{i\neq j}W(x_i-x_j)
\]
is usually called energy or Hamiltonian of the system. We assume that $\be_N$,
$V$, and $W$ are chosen in such a way that $Z_N<\infty$ for any $N$. The law
$P_N$ is invariant by permutation of the coordinates $x_1,\ldots,x_N$
(exchangeable), and $H_N$ depends only on the empirical measure
\[
  \mu_N=\frac{1}{N}\sum_{i=1}^N\de_{x_i}.
\]
Therefore $P_N$ is also the law of a random empirical measure encoding a cloud
of indistinguishable particles $x_1,\hdots,x_N$. We emphasize that the
particles live on the space $S^N=S\times\cdots\times S$ of dimension $dN$. The
parameter $n$ serves as the physical dimension of the ambient space, for the
Coulomb gas setting described next.

For any $m\geq1$ and $x\in\dR^m$, we denote by
$\ABS{x}=\sqrt{x_1^2+\cdots+x_m^2}$ the Euclidean norm of $x$. This matches
the absolute value when $m=1$ and the modulus when $m=2$, $\dR^2\equiv\dC$.

\subsection{Coulomb gases}

The notion of Coulomb gas is based on elementary electrostatics. Here the
vector subspace $S$ is interpreted as a conductor. It corresponds to taking
$W=g$ where $g$ is the Coulomb kernel or Green function in the physical space
$\dR^n$. More precisely, recall that the Green function $g$ in $\mathbb{R}^n$,
$n\geq2$, is defined for all $x\in\mathbb{R}^n$, $x\neq0$, by
\[
  g(x)=
  \begin{cases}
    \log\frac{1}{\ABS{x}} & \text{if $n=2$,}\\
    \frac{1}{\ABS{x}^{n-2}} & \text{if $n\geq3$}.
  \end{cases}
\]
This function is the fundamental solution of the Poisson equation, namely,
denoting by $\De$ the Laplace operator in $\dR^n$ and by $\de_0$ the Dirac mass at
$0$, we have, in the sense of distributions,
\[
  -\Delta g=c\delta_0,
\quad\mathrm{with}\quad
c
=\begin{cases}
  2\pi & \text{if $n=2$,}\\
  (n-2)|\mathbb{S}^{n-1}|=\frac{n(n-2)\pi^{n/2}}{\Gamma(1+n/2)} & \text{if $n\geq3$}.    
\end{cases}
\]
The physical interpretation in terms of electrostatics is as follows:
$H_N(x_1,\ldots,x_N)$ is the electrostatic energy of a configuration of $N$
electrons in $\dR^n$ lying on $S$ at positions $x_1,\ldots,x_N$, in an
external field given by the potential $V$. The Green function or Coulomb
kernel $g$ expresses the Coulomb repulsion which is a two body singular
interaction. The probability measure $P_N$ can be seen as a Boltzmann--Gibbs
measure, $\be_N$ playing the role of an inverse temperature. The probability
measure $P_N$ is known as a Coulomb gas or as a one-component plasma, see for
instance \cite{serfaty-icm2018} and references therein.

\subsection{Log-gases} A log-gas corresponds to choosing $d=n$ and a
logarithmic interaction $W$ whatever the value of $n$ is, namely
\[
  W(x)=\log\frac{1}{\ABS{x}}=-\frac{1}{2}\log(x_1^2+\cdots+x_d^2),\quad x\in S.
\]
Coulomb gases and log-gases coincide when $d=n=2$. In dimension $d=n\geq3$,
log-gases are natural and classical objects of approximation theory and can be
seen as limiting Riesz potentials, namely
$\lim_{\alpha\to0}\frac{1}{\alpha}(|x|^{-\alpha}-1)$, see for instance
\cite{MR1754783,MR1631413,serfaty-icm2018}.

\subsection{Static energy and equilibrium measures}

Under natural assumptions over $V$ and $W$, typically when $\be_N\gg N$ and
$V$ beats $W$ at infinity, it is well known, see for
instance~\cite{MR3262506,MR3309890} and references therein, that $P_N$ almost surely,
the empirical measure
\[
  \mu_N=\frac{1}{N}\sum_{i=1}^N\de_{x_i}
\]
tends as $N\to\infty$ to a non random probability measure, the equilibrium
measure
\[
  \mu_*=\arg\inf\cE,
\]
the unique minimizer of the strictly convex lower semi-continuous ``energy''
$\cE$ defined by
\[
  \mu\mapsto\cE(\mu)=\int\!V\dd\mu+\iint W(x-y)\mu(\dd x)\mu(\dd y).
\]
When $W=g$ is the Coulomb kernel, the quantity $\cE(\mu)$ is the electrostatic
energy of the distribution of charges $\mu$, formed by the sum of the
electrostatic potential coming from the external electric field $V$ with the
Coulomb self repulsion by mean of the Coulomb kernel $g$. Note that
$\cE(\mu)=\infty$ if $\mu$ has a Dirac mass due to the singularity of $g$. An
Euler--Lagrange variational analysis reveals that when $S=\dR^d$ and $V$ is
smooth, convex, and grows faster than $g$ at infinity then the equilibrium
probability measure $\mu_*$ is compactly supported and has density
proportional to $\Delta V$, see \cite{MR3262506} and references therein.
Table~\ref{fi:equilibrium} gives examples of equilibrium measures in this
Coulomb setting. We refer to
\cite{MR0350027,MR3308615,MR1485778,MR3309890,serfaty-icm2018} for old and new
potential theory from this analytic point of view. Moreover, quite a few equilibrium
measures are known for log-gases beyond Coulomb gases, see for instance~\cite{chafai-saff}.

\begin{table}[h]
  \[
    \begin{array}{c|c|c|c|c|c}
      d & S & n & V & \mu_* & \text{Nickname}\\\hline\hline%
      1 & \dR & 2 & \infty\mathbf{1}_{\text{interval}^c} & \text{arcsine} & \\%
      1 & \dR & 2 & x^2 & \text{semicircle} & \text{GUE}\\%
      2 & \dR^2 & 2 & |x|^2 & \text{uniform on a disc} & \text{Ginibre}\\%
      d\geq3 & \dR^d & d & | x|^2 & \text{uniform on a ball} & \\%
      d\geq3 & \dR^d & d & \text{radial} & \text{radial in a ring} & 
    \end{array}
  \]
  \caption{\label{fi:equilibrium} Examples of equilibrium measures for Coulomb
    gases, see \cite{MR1485778,MR3262506}.}
\end{table}

Actually it can be shown that essentially if $\be_N\gg N$ and $V$ beats $g$ at
infinity then under ${(P_N)}_N$ the sequence of random empirical measures
${(\mu_N)}_N$ satisfies a large deviation principle with speed $\be_N$ and
good rate function $\cE$, see \cite{MR3262506,garcia-zelada,berman}.
Concentration of measure inequalities are also available, see \cite{coco} and
references therein.

\subsection{Two remarkable gases from random matrix theory}
\label{sec:remarkable}
Let us give a couple of famous gases from random matrix theory that will serve
as benchmark for our algorithm. They correspond to $n=2$ because the Lebesgue
measure on a matrix translates via the Jacobian of the change of variable to a
Vandermonde determinant on the eigenvalues, giving rise to the two-dimensional
Coulomb kernel inside the exponential via the identity
\[
  \prod_{i<j}|x_i-x_j|=\exp\Bigr(\sum_{i<j}\log|x_i-x_j|\Bigr).
\]
Hence the name ``log-gases''. A good reference on this subject
is~\cite{MR2641363} and we refer
to~\cite{MR2168344,MR1936554,MR2641363,MR3458536,MR3615091} for more examples
of Coulomb gases related to random matrix models. Coulomb gases remain
interesting in any dimension~$n$ beyond random matrices,
see~\cite{MR3309890,serfaty-icm2018}.

\textbf{Beta-Hermite model.} This model corresponds to 
\[
  d=1,
  \ n=2,
  \ S=\mathbb{R},
  \ V(x)=\frac{x^2}{2\be},
  \ W(x)=-\log\ABS{\cdot},
  \ \be_N=N^2\be,
  \ \be\in(0,\infty).
\]
This means that the particles evolve on the line $\dR$ with Coulomb
interactions given by the Coulomb kernel in $\dR^2$. For $\be=2$, it becomes
the famous Gaussian Unitary Ensemble (GUE), which is the distribution of the
eigenvalues of random $N\times N$ Hermitian matrices distributed according to
the Gaussian probability measure with density proportional to
$H\mapsto\e^{-N\mathrm{Tr}(H^2)}$. Beyond the case $\be=2$, the cases $\be=1$
and $\be=4$ correspond respectively to Gaussian random matrices with real and
quaternionic entries. Following~\cite{MR1936554}, for all $\be\in(0,\infty)$,
the measure~$P_N$ is also the distribution of the eigenvalues of special
random $N\times N$ Hermitian tridiagonal matrices with independent but non
identically distributed entries. Back to the case~$\be=2$, the law~$P_N$
writes
\begin{equation}\label{eq:BH}
 (x_1,\ldots,x_N)\in\dR^N\mapsto\e^{-\frac{N}{2}\sum_{i=1}^N x_i^2}\prod_{i<j}(x_i-x_j)^2.
\end{equation}
In this case, the Coulomb gas $P_N$ has a determinantal structure, making it
integrable or exactly solvable for any $N \geq 2$,
see~\cite{MR2129906,MR2641363}. This provides in particular a formula for the
density of the mean empirical spectral distribution $\dE\mu_N$ under $P_N$,
namely
\begin{equation}\label{eq:EGUEN}
  x\in\dR\mapsto
  \frac{\e^{-\frac{N}{2}x^2}}{\sqrt{2\pi N}}
  \sum_{\ell=0}^{N-1}H_\ell^2(\sqrt{N}x),
\end{equation}
where ${(H_\ell)}_{\ell\geq0}$ are the Hermite polynomials which are the
orthonormal polynomials for the standard Gaussian distribution $\cN(0,1)$. The
equilibrium measure $\mu_*$ in this case is the Wigner semicircle distribution
with the following density with respect to the Lebesgue measure:
\begin{equation}\label{eq:SC}
  x\in\dR\mapsto\frac{\sqrt{4-x^2}}{2\pi}\mathbf{1}_{x\in[-2,2]}.
\end{equation}
A plot of $\mu_*$ and $\dE\mu_N$ is provided in
Figure~\ref{fi:GUE}, together with our simulations. We refer
to~\cite{MR2041832} for a direct proof of convergence of~\eqref{eq:EGUEN}
to~\eqref{eq:SC} as $N\to\infty$. Beyond the case $\be=2$, the equilibrium
measure $\mu_*$ is still a Wigner semicircle distribution, scaled by $\be$,
supported by the interval $[-\be,\be]$, but up to our knowledge we do not have
a formula for the mean empirical spectral distribution $\dE\mu_N$, except when
$\be$ is an even integer, see~\cite{MR1936554}.

\textbf{Beta-Ginibre model.} This model corresponds to
\[
  d=2,
  \ n=2,
  \ S=\mathbb{R}^2,
  \ V(x)=\frac{|x|^2}{\be},
  \ W(x)=-\log\ABS{x},
  \ \be_N=N^2\be,
  \ \be\in(0,\infty).
\]
In this case, the particles move in $\dR^2$ with a Coulomb repulsion of
dimension $2$ -- it is therefore a Coulomb gas. As for the GUE, the law $P_N$
can be written as
\begin{equation}\label{eq:BG}
  (x_1,\ldots,x_N)\in(\dR^2)^N
  \mapsto\e^{-N\sum_{i=1}^N|x_i|^2}\prod_{i<j}|x_i-x_j|^\be.
\end{equation}

When~$\be=m$ for an even integer $m\in\{2,4,\ldots\}$, the law of this gas
matches the Laughlin wavefunction modeling the fractional quantum Hall effect
(FQHE), see for instance \cite{MR2377026}.

For $\be=2$, this gas, known as the complex Ginibre Ensemble, matches the
distribution of the eigenvalues of random $N\times N$ complex matrices
distributed according to the Gaussian probability measure with density
proportional to $M\mapsto\e^{-N\mathrm{Tr}(MM^*)}$ where
$M^*=\overline{M}^\top$. In this case $P_N$ has a determinantal structure,
see~\cite{MR2129906,MR2641363}. This provides a formula for the density of the
mean empirical spectral distribution $\dE\mu_N$ under $P_N$, namely
\begin{equation}\label{eq:EGN}
  x\in\dR^2\mapsto
  \frac{\e^{-N|x|^2}}{\pi}\sum_{\ell=0}^{N-1}\frac{|\sqrt{N}x|^{2\ell}}{\ell!},
\end{equation}
which is the analogue of~\eqref{eq:EGUEN} for the Gaussian Unitary Ensemble.
Moreover, if $Y_1,\ldots,Y_N$ are independent and identically distributed
Poisson random variables of mean $|x|^2$ for some $x\in\dR^2$,
then~\eqref{eq:EGN} writes
\[
  x\in\dR^2\mapsto\frac{1}{\pi}\dP\left(\frac{Y_1+\cdots+Y_N}{N}<1\right).
\]
As $N\to\infty$, by the law of large numbers, it converges to $1/\pi$ if
$|x|<1$ and to $0$ if $|x|>1$, while by the central limit theorem it converges
to $1/(2\pi)$ if $|x|=1$. It follows that $\dE\mu_N$ converges weakly as
$N\to\infty$ to the uniform distribution on the disk, with density
\begin{equation}\label{eq:G}
  x\in\dR^2 \mapsto \frac{\mathbf{1}_{|x|<1}}{\pi},
\end{equation}
which is the equilibrium measure $\mu_*$. When $N$ is finite, the numerical
evaluation of~\eqref{eq:EGN} is better done by mean of the Gamma law. Namely,
by induction and integration by parts,~\eqref{eq:EGN} writes
\[
  x\in\dR^2 \mapsto
  \frac{1}{\pi(N-1)!}\int_{N|x|^2}^\infty u^{N-1}\e^{-u}\dd u
  =\frac{\Gamma(N,N|x|^2)}{\pi},
\]
where $\Gamma$ is the normalized incomplete Gamma function and where we used
the identity
\[
  \e^{-r}\sum_{\ell=0}^{N-1}\frac{r^\ell}{\ell!}
  =
  \frac{1}{(N-1)!}\int_r^\infty u^{N-1}\e^{-u}\dd u.
\]
Note that $t\mapsto 1-\Gamma(N,t)$ is the cumulative distribution function of
the Gamma distribution with shape parameter $N$ and scale parameter $1$.
Figure~\ref{fi:Ginibre} illustrates the
difference between the limiting distribution~\eqref{eq:G} and the mean
empirical spectral distribution~\eqref{eq:EGN} for a finite $N$. Beyond the
case $\be=2$, we no longer have a formula for the density of $\dE\mu_N$, but a
simple scaling argument reveals that the equilibrium measure $\mu_*$ is in
this case the uniform distribution on the centered disk of radius
$\sqrt{\beta/2}$.

\section{Simulating log-gases and Coulomb gases}

Regarding simulation of log-gases or Coulomb gases such as~\eqref{eq:gibbs},
it is natural to use the random matrix models when they are available. There
exist also methods specific to determinantal processes which cover the
log-gases of random matrix theory with~$\be=2$,
see~\cite{MR2216966,MR2551206,MR3334666,MR3439168,bardenet-hardy,MR3382600,hardy}.
Beyond these specially structured cases, a great variety of methods are
available for simulating Boltzmann--Gibbs measures, such as overdamped
Langevin diffusion algorithm, Metropolis--Hastings algorithm, Metropolis
adjusted Langevin algorithm (MALA), and kinetic versions called Hybrid or
Hamiltonian Monte Carlo (HMC) which are based on a kinetic (or underdamped)
Langevin diffusion, see for instance~\cite{MR2742422,MR3509213}. Other
possibilities exist, such as Nos\'{e}-Hoover dynamics~\cite{jones2011adaptive}
or piecewise deterministic Markov processes~\cite{bouchard2018bouncy}.

Two difficulties arise when sampling measures as~\eqref{eq:gibbs}. First, the
Hamiltonian $H_N$ involves all couples, so the computation of forces and
energy scales quadratically with the number of particles. A natural way to
circumvent this numerical problem is to use clusterization procedures such as
the ``fast multipole methods'', see for instance~\cite{MR3659634}. A second
difficult feature of such a Hamiltonian is the singularity of the interacting
function $W$, which typically results in numerical instability. A standard
stabilization procedure is to <<tame>> the
dynamics~\cite{hutzenthaler2012strong,brosse2017tamed}, which is the strategy
adopted in~\cite{MR3124976}. However, this smoothing of the force induces a
supplementary bias in the invariant measure, as shown
in~\cite{brosse2017tamed} for regular Hamiltonians. This requires using small
time steps, hence long computations. In the present note, we explore for the
first time the usage of HMC for general Coulomb gases in the context of random
matrices, in the spirit of~\cite{stoltz-trstanova}, the difficulty being the
singularity of the interaction. This method has the advantage of sampling the
exact invariant measure~\eqref{eq:gibbs}, while allowing to choose large time
steps, which reduces the overall computational cost~\cite{fathi2015error}.

In Section~\ref{sec:review}, we review standard methods for sampling measures
of the form~$\e^{-\beta_N H_N}$, before presenting in detail the HMC algorithm
in Section~\ref{sec:HMC}.

\subsection{Standard sampling methods}
\label{sec:review}

To simplify and from now on, we suppose the support set $S$
in~\eqref{eq:gibbs} to be~$\dR^d$. We introduce the methods based on the
overdamped Langevin dynamics. To sample approximately~\eqref{eq:gibbs}, the
idea is to exploit the fact that $P_N$ in~\eqref{eq:gibbs} is the reversible
invariant probability measure of the Markov diffusion process
${(X_t)}_{t\geq0}$ solution to the stochastic differential equation:
\begin{equation}\label{eq:overdamped}
  \dd X_t=-\al_N\nabla H_N(X_t)\,\dd t+\sqrt{2\frac{\al_N}{\be_N}}\,\dd B_t,
\end{equation}
or in other words
\[
  X_t=X_0-\al_N\int_0^t\na H_N(X_s)\,\dd s+\sqrt{2\frac{\al_N}{\be_N}}B_t,
\]
where ${(B_t)}_{t\geq0}$ is a standard Brownian motion on $S^N$ and $\al_N>0$
is an arbitrary time scaling parameter (for instance $\al_N=1$ or
$\al_N=\be_N$). The infinitesimal generator associated
with~\eqref{eq:overdamped} is
\[
  Lf=\frac{\al_N}{\be_N}\De f -\al_N\na H_N\cdot\na f.
\]
The difficulty in solving \eqref{eq:overdamped} lies in the fact that the
energy $H_N$ involves a singular interaction $W$, which may lead the process
to explode. Actually, under certain conditions on $\be_N$ and $V$, the
equation \eqref{eq:overdamped} is well posed, the process ${(X_t)}_{t\geq0}$
is well defined, and
\[
X_t \,
\underset{t\to\infty}{\overset{\mathrm{Law}}{\longrightarrow}} \,
P_N,
\]
for all non-degenerate initial condition $X_0$. See for instance
\cite{MR2760897,MR3699468,chafai-lehec} for the case of Beta-Hermite case
known as the Dyson Ornstein--Uhlenbeck process, and~\cite{coulsim} for the
Beta-Ginibre case. We do not discuss these delicate aspects in this note. A
convergence in Cesàro mean is provided by the ergodic theorem for additive
functionals,
\[
\frac{1}{t}\int_0^t\de_{X_s}\, \dd s\,
\underset{t\to\infty}{\overset{\mathrm{weak}}{\longrightarrow}} \,
P_N
\]
almost surely or, for any test function $f\in L^1(P_N)$,
\[
\frac{1}{t} \int_0^tf(X_s)\,\dd s\, \underset{t\to\infty}{\longrightarrow}\,
\int_S f \, \dd P_N,
\]
almost surely. It is also possible to accelerate the convergence by adding a
divergence free term in the dynamics~\eqref{eq:overdamped}, see for instance
\cite{2016JSP...163..457D,2013JSP...152..237L} and references therein. This
modification keeps the same invariant distribution but produces a
non-reversible dynamics.

This method of simulation is referred to as an ``unadjusted Langevin
algorithm'', a terminology which will be clarified later on. In practice, one
cannot simulate the continuous stochastic process ${(X_t)}_{t\geq0}$ solution
to~\eqref{eq:overdamped}, and resorts to a numerical integration with a finite
time step $\Dt$. A typical choice is the Euler--Maruyama
scheme~\cite{MR1260431,milstein2013stochastic}, which reads
\begin{equation}\label{eq:euler}
  x_{k+1} = x_k -\na H_N(x_k)\al_N\Dt + \sqrt{2\frac{\al_N}{\be_N}\Dt}G_k,  
\end{equation}
where $(G_k)$ is a family of independent and identically distributed standard
Gaussian variables, and $x_k$ is an approximation of $X_{k\Dt}$. Note that
$\al_N$ and $\Dt$ play the same role here. However, because of the singularity
of $H_N$, this sampling scheme leads to important biases in practice,
and~\eqref{eq:euler} may even lack an invariant
measure~\cite[Section~6]{MATTINGLY2002185}. One way to stabilize the dynamics
is to use a tamed version of~\eqref{eq:euler}, which typically takes the
following form:
\begin{equation}
  \label{eq:tamed}
  x_{k+1} = x_k %
  - \frac{\na H_N(x_k)\al_N\Dt}
  {1+|\na H_N(x_k)|\al_N\Dt} + \sqrt{2\frac{\al_N}{\be_N}\Dt}G_k.
\end{equation}
This strategy is used in~\cite{MR3124976} but, as noted by the authors, the
number of time steps needed to run a trajectory of fixed time $T$ scales as
$\Dt\sim N^{-2}$, which makes the study of large systems difficult.

Another strategy is to add a selection step at each iteration. This is the
idea of the Metropolis Adjusted (overdamped) Langevin Algorithm
(MALA)~\cite{roberts1996exponential}, which
prevents irrelevant moves with a Metropolis step. One can also view the MALA
algorithm as a Metropolis algorithm in which the proposal is produces by using
a one step discretization of the Langevin dynamics~\eqref{eq:overdamped}.
Let us make this precise; more details can be found \textit{e.g.}
in~\cite{roberts1996exponential,MR2080278}.

\begin{algo}[Metropolis Adjusted (overdamped) Langevin Algorithm --  MALA]
  \label{algo:MALA}
  Let $K$ be the Gaussian transition kernel associated to the Markov chain of
  the Euler discretization~\eqref{eq:euler} of the
  dynamics~\eqref{eq:overdamped}. For each step $k$,
  \begin{itemize}
  \item draw a proposal $\tilde{x}_{k+1}$ according to the kernel
    $K(x_k,\cdot)$,
  \item compute the probability 
    \begin{equation}\label{eq:ratio}
      p_k=1\wedge
      \frac{K(\tilde{x}_{k+1},x_k)\e^{-\be_N H_N(\tilde{x}_{k+1})}}
      {K(x_{k},\tilde{x}_{k+1})\e^{-\be_N H_N(x_{k})} },
    \end{equation}
  \item set
    \[
      x_{k+1}
      =\begin{cases}
        \tilde{x}_{k+1} & \text{with probability $p_k$};\\
        x_k & \text{with probability $1-p_k$}.
      \end{cases}
    \]
  \end{itemize}  
\end{algo}

Note that the ``reversed'' kernel $K(\cdot,x)$ is Gaussian only if~$H_N$ is a
quadratic form. Note also that if the proposal kernel~$K$ is symmetric in the
sense that $K(x,y)=K(y,x)$ for all $x,y$ then it disappears
in~\eqref{eq:ratio}, and it turns out that this is the case for the Hybrid
Monte Carlo algorithm described next (up to momentum reversal)!

A natural issue with these algorithms is the choice of $\Dt$: if it is too
large, an important fraction of the proposed moves will be rejected, hence
poor convergence properties; conversely, if $\Dt$ is too small, many steps
will be accepted but the physical ellapsed time will be small, hence a large
variance for a fixed number of iterations. This algorithm actually has a nice
scaling of the optimal time step $\Dt$ with the dimension of the system.
Indeed, it can be shown that it scales as $\Dt \sim N^{-\frac{1}{3}}$, at
least for product measures (see~\cite{roberts2001optimal} and references
therein). Although this algorithm is already efficient, we propose to use a
kinetic version with further advantages.

\subsection{Hybrid Monte Carlo algorithm}
\label{sec:HMC}

Hybrid Monte Carlo is built on Algorithm~\ref{algo:MALA}, but using a kinetic
version of~\eqref{eq:overdamped}. For this, a momentum variable is introduced so
as to improve the exploration of the space. Namely, set $E=\dR^{dN}$, and let
$U_N:E\to\dR$ be smooth and such that $\e^{-\be_N U_N}$ is Lebesgue
integrable. Let ${(X_t,Y_t)}_{t\geq0}$ be the diffusion process on $E\times E$
solution to the stochastic differential equation
\begin{equation}\label{eq:underdamped}
  \setlength{\jot}{5pt}
  \left\{
  \begin{aligned}
    \dd X_t & =\al_N\na U_N(Y_t)\,\dd t,\\ 
    \dd Y_t &\displaystyle =-\al_N\na H_N(X_t)\,\dd t-\ga_N\al_N\na U_N(Y_t)\,\dd
    t+\sqrt{2\frac{\ga_N\al_N}{\be_N}}\,\dd B_t,
  \end{aligned}
  \right.
\end{equation}
where ${(B_t)}_{t\geq0}$ is a standard Brownian motion on $E$, and $\ga_N>0$
is an arbitrary parameter which plays the role of a friction, and which may
depend a priori on $N$ and $(X_t)_{t\geq 0}$, even if we do not use this
possibility here. In addition, $H_N$ and $\be_N$ are as in \eqref{eq:gibbs},
while $U_N$ plays the role of a generalized kinetic
energy~\cite{stoltz-trstanova}. This dynamics admits the following generator:
\begin{equation}\label{eq:L1L2}
  Lf=
  \underbrace{-\al_N\na H_N(x)\cdot\na_yf+\al_N\na U_N(y)\cdot\na_x f}_{L_1}
  +
  \underbrace{\frac{\ga_N\al_N}{\be_N}\De_yf-\ga_N\al_N\na U_N(y)\cdot\na_yf}_{L_2}
\end{equation}
where $L_1$ is known as the Hamiltonian part while $L_2$ is called the
fluctuation-dissipation part. The dynamics leaves invariant the product
Boltzmann--Gibbs measure
\[
  R_N=P_N\otimes Q_N
  \quad\text{where}\quad
  Q_N(\dd y)=\frac{\e^{-\be_NU_N(y)}}{Z_N'}\dd y,
\]
see for instance \cite{stoltz-trstanova}. In other words
\begin{equation}\label{eq:equilangevin}
  R_N(\dd x,\dd y)=\frac{\e^{-\be_N\widetilde H_N(x,y)}}{Z_NZ_N'}\dd x\, \dd y
  \quad\text{with}\quad
  \widetilde H_N(x,y)=H_N(x)+U_N(y).
\end{equation}
As for the overdamped dynamics, the ergodic theorem for additive functionals gives
\[
  \frac{1}{t}\int_0^t\de_{(X_s,Y_s)}\, \dd s\,
  \underset{t\to\infty}{\overset{\mathrm{weak}}{\longrightarrow}}\, R_N
  \quad\text{almost surely.}
\]

\begin{remark}[Terms: Hamiltonian, Langevin, overdamped, underdamped, kinetic]
  The dynamics~\eqref{eq:underdamped} is called ``Hamiltonian'' when we turn
  off the noise by taking $\ga_N=0$. On the other hand, when $\ga_N\to\infty$
  and $\al_N\to0$ with $\al_N\ga_N=1$, we recover \eqref{eq:overdamped} from
  \eqref{eq:underdamped} with $Y_t$ and $U_N$ instead of $X_t$ and $H_N$. Both
  \eqref{eq:overdamped} and \eqref{eq:underdamped} are known as Langevin
  dynamics. To be more precise,~\eqref{eq:overdamped} is generally called
  overdamped while~\eqref{eq:underdamped} is referred to as kinetic or
  underdamped.
\end{remark}

When $U_N(y)=\frac{1}{2}|y|^2$ then $Y_t=\dd X_t/\dd t$, and in this case
$X_t$ and $Y_t$ can be interpreted respectively as the \emph{position} and the
\emph{velocity} of a system of $N$ points in $S$ at time $t$. In this case we
say that $U_N$ is the \emph{kinetic energy}. For simplicity, we specialize in
what follows to this ``physical'' or ``kinetic'' case and refer
to~\cite{stoltz-trstanova} for more possibilities.

As before, to simulate ${(X_t,Y_t)}_{t\geq0}$, one can
discretize~\eqref{eq:underdamped} and sample from a trajectory. This will
provide a proposal for the HMC scheme as the Euler
discretization~\eqref{eq:euler} did for Algorithm~\ref{algo:MALA}. A good way
of doing this is a splitting procedure. First, one integrates the Hamiltonian
part \textit{i.e.} the operator $L_1$ in~\eqref{eq:L1L2}, which amounts to a
standard Hamiltonian dynamics, before integrating the fluctuation-dissipation
part \textit{i.e.} the operator $L_2$ in~\eqref{eq:L1L2}. For discretizing the
Hamiltonian dynamics over a time step, a standard approach is the Verlet
integrator~\cite{hairer2006geometric,MR2681239}, which we describe now. For a
time step $\Dt >0$, this scheme reads, starting from a state $(x_k,y_k)$ at
time $k$:
\[  
\left\{
\begin{aligned}
  y_{k+\frac{1}{2}} &=  y_k -\na H_N(x_k)\al_N\frac{\Dt}{2}, \\
  x_{k+1} &= x_k +y_{k+\frac{1}{2}}\al_N\Dt, \\
  \tilde{y}_{k+1} &= y_{k+\frac{1}{2}}-\na H_N(x_{k+1})\al_N\frac{\Dt}{2}.
\end{aligned}
\right.
\]
This corresponds to updating the velocity over half a time step, then the
positions over a time step, and again the velocity over half a time-step.
Given that this scheme only corresponds to the Hamiltonian part, it remains to
integrate the fluctuation-dissipation part, corresponding to $L_2$
in~\eqref{eq:L1L2}. For quadratic energies, it is a simple Ornstein--Uhlenbeck
process whose variance can be computed explicitly. Therefore, we add to the
previous scheme the following velocity update which comes from the Mehler
formula\footnote{The Mehler formula states that the Ornstein--Uhlenbeck
  process ${(Z_t)}_{t\geq0}$ in $\dR^n$ solution of the stochastic
  differential equation %
  $\dd Z_t=\sqrt{2\sigma^2}\dd B_t-\rho Z_t\dd t$ satisfies
  $\mathrm{Law}(Z_{t+s}\mid Z_s=z)=\mathcal{N}(z\e^{-\rho
    t},\frac{1-\e^{-2\rho t}}{\rho}\sigma^2I_n)$.}:  
\[
y_{k+1}= \eta\tilde{y}_{k+1} +\sqrt{\frac{1-\eta^2}{\be_N}}G_k,
\quad
\eta=\e^{-\ga_N\al_N\Dt},
\]
where $G_k$ is a standard Gaussian random variable. Like the numerical
scheme~\eqref{eq:euler}, because of the singularity of the interactions, this
integrator may not have an invariant measure~\cite{MATTINGLY2002185}, or its
invariant measure may be a poor approximation of $R_N$ depending on the time
step~\cite{leimkuhler2015computation}. Note that, here again, $\al_N$ and
$\Dt$ play the same role.

Hybrid or Hamiltonian Monte Carlo (HMC) methods, built on the later
integration, appeared in theoretical physics in lattice quantum chromodynamics
with~\cite{duane-hmc}, see also \cite{rossky-doll-friedman}, and are still
actively studied in applied mathematics, see for
instance~\cite{MR3129023,stoltz-trstanova,MR2681239,MR3214779,bou2017geometric,durmus2017convergence,dalalyan-riou-durand}
and references therein. The HMC algorithm can be thought of in a sense as a
special Metropolis Adjusted (underdamped) Langevin Algorithm. Indeed, inspired
by the MALA Algorithm~\ref{algo:MALA}, a way to avoid the stability problem of
the discretization of the kinetic Langevin dynamics mentioned above is to add
an acceptance-rejection step. A surprising advantage of this approach is that
the Verlet integration scheme is time reversible up to momenta
reversal~\cite[Sec.\ 2.1.3 and eq.\ (2.11)]{MR2681239}, hence when computing
the acceptance probability as in~\eqref{eq:ratio}, the transition kernel does
not appear. Note that the Verlet algorithm has been widely used for years by
statistical physicists, and goes back to the historical works of Verlet
\cite{verlet} and Levesque and Verlet
\cite{levesque-verlet-bis,levesque-verlet}. Let us now describe the algorithm.

\begin{algo}[HMC]\label{algo:HMC}
  Start from a configuration $(x_0,y_0)$ and perform the following steps for
  each time $k\geq 0$:
  \begin{enumerate}
  \item update the velocities with
    \[
      \tilde{y}_{k}=\eta y_{k}+\sqrt{\frac{1-\eta^2}{\beta_N}}
      G_k, \quad \eta = \e^{-\ga_N\al_N\Dt};
    \]
  \item run one step of the Verlet scheme:
    \begin{equation}\label{eq:verlet}
      \left\{
      \begin{aligned}
        \tilde{y}_{k+\frac{1}{2}}
        &=\tilde{y}_k -\na H_N(x_k)\al_N\frac{\Dt}{2};\\
        \tilde{x}_{k+1}
        &=x_k +\tilde{y}_{k+\frac{1}{2}}\al_N\Dt;\\
        \tilde{y}_{k+1}
        &=\tilde{y}_{k+\frac{1}{2}}-\na H_N(x_{k+1})\al_N\frac{\Dt}{2};\\
      \end{aligned}
    \right.
    \end{equation}
  \item compute the probability ratio
    \[
      p_k = 1\wedge\exp\left[-\be_N\Big(H_N(\tilde{x}_{k+1}) +
          \frac{\tilde{y}_{k+1}^2}{2} - H_N(x_{k}) - \frac{\tilde{y}_{k}^2}{2}
        \Big)\right];
    \]
  \item set
    \[
      (x_{k+1},y_{k+1})
      =
      \begin{cases}
        (\tilde{x}_{k+1},\tilde{y}_{k+1})
        & \text{with probability $p_k$};\\
        (x_{k},-\tilde{y}_{k})
        & \text{with probability $1-p_k$.}
      \end{cases}
    \]
  \end{enumerate}
\end{algo}

As noted in the various references above, the Metropolis step acts as a
corrector on the energy conservation of the Hamiltonian step. In this, it
helps avoiding irrelevant moves, while enhancing the exploration capacities of
the dynamics through the speed variable. A more precise argument in favor of
this algorithm is the scaling of the time step $\Dt$ with respect to the
system size $N$. Indeed, as shown in~\cite{MR3129023} for product measures,
the optimal scaling is as $\Dt\sim N^{-\frac{1}{4}}$, which makes the
algorithm appealing for large systems. Since the Hamiltonian computational
cost scales as $N^2$, we see that the cost of the algorithm for a fixed time $T$
and $N=\ceil{T/\Dt}$ is in $\mathcal{O}(N^\frac{9}{4})$, which has to be
compared to the $\mathcal{O}(N^4)$ cost reached in~\cite{MR3124976}. Finally, the
parameter $\ga_N$ can also be tuned in order to optimize the speed of
convergence -- we leave this point here and stick to $\ga_N=1$.

The control of the error or rate of convergence for the HMC algorithm is the
subject of active research, see for instance~\cite{lee-vempala}
and~\cite{durmus2017convergence,bourabee-eberle-zimmer} for some results under
structural assumptions.

From a practical point of view, the algorithm can be tested in the following
way. First, when only the Hamiltonian part of the dynamics is integrated with
the Verlet scheme \eqref{eq:verlet}, it can be checked that the energy
variation over one time step scales as $\Dt^3$ as $\Dt\to 0$. Then, if the
selection step is added, the rejection rate should also scale as $\Dt^3$. When
the momentum resampling is added, this rejection rate scaling should not
change. For completeness, we illustrate some of these facts in
Section~\ref{sec:numerics}.

\section{Numerical experiments on remarkable models}
\label{sec:numerics}

In this section, we start testing Algorithm~\ref{algo:HMC} for the two cases
described in Section~\ref{sec:remarkable}. Since the equilibrium measures are
known for any $N\geq 2$, we will be able to compare accurately our results
with the expected one. We will also consider models for which the empirical
spectral distribution and the equilibrium distribution are not known. We
remind that when $S=\dR^d$ with $d\geq 1$ we have the following formulas that
hold in any dimension:
\[
  \nabla \ABS{x}^2=2x,\quad
  \nabla\log\frac{1}{\ABS{x}}
  =-\dfrac{x}{\ABS{x}^2},\quad
  \nabla\frac{1}{\ABS{x}}=-\dfrac{x}{\ABS{x}^3}.
\]

\subsection{Case study: 1D}
\label{sec:GUE}

\begin{figure}[p]
  \centering
  \includegraphics[width=.8\textwidth]{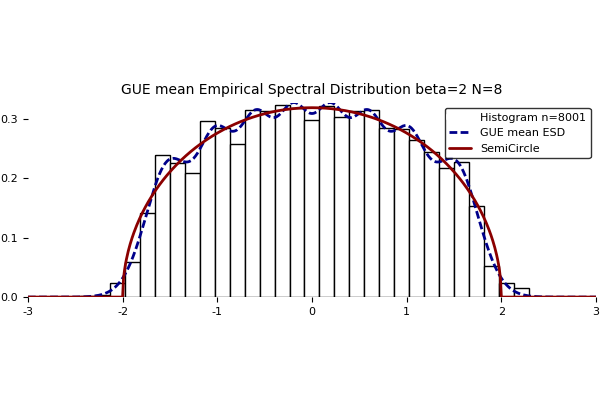}\\
  \includegraphics[width=.8\textwidth]{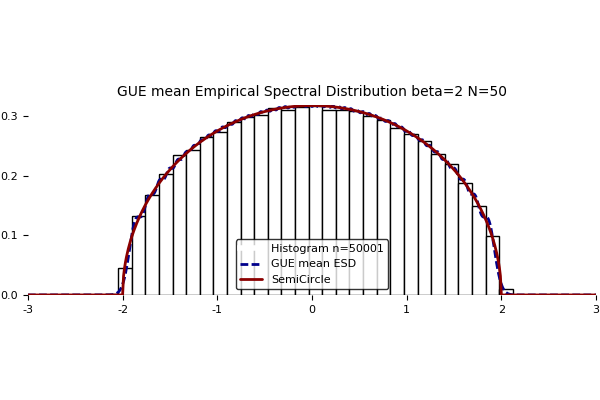}  
  \caption{\label{fi:GUE} Study of the Gaussian Unitary Ensemble with $N=8$
    (top) and $N=50$ (bottom). The solid line is the plot of the limiting
    spectral distribution~\eqref{eq:SC} while the dashed line is the plot of
    the mean empirical distribution~\eqref{eq:EGUEN}. The bars form the
    histogram of simulations obtained using our HMC algorithm. This algorithm
    was run once with final-time $T=10^6$ and time-step $\Dt=0.5$. The
    histogram was produced by looking at the last half of the~trajectory and
    retaining the positions each~$1000$ time-steps, producing $n$ values,
g    namely $\approx~8\times 10^3$ and $\approx~5\times 10^4$ respectively.}
\end{figure}

\begin{figure}[p]
  \centering
  \includegraphics[width=.8\textwidth]{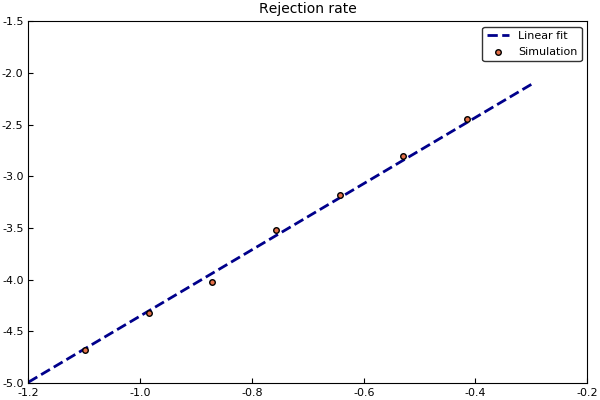}\\
  \caption{\label{fi:rejection} Evolution of the rejection rate in
    Algorithm~\ref{algo:HMC} as~$\Dt$ goes to zero, for the Gaussian Unitary
    Ensemble with~$N=50$, $\beta=2$ and $T=10^{5}$ (in log-log coordinate). }
\end{figure}

We test the numerical method by looking at the mean empirical distribution in
the case of the Gaussian Unitary Ensemble~\eqref{eq:BH} with $\be=2$, $N=8$,
for which the exact expression of $\dE\mu_N$ under $P_N$ is provided
by~\eqref{eq:EGUEN}. The results in Figure~\ref{fi:GUE} show a very good
agreement between the exact result and the algorithm. For completeness, we
study the rejection rate of the algorithm as~$\Dt$ goes to zero, as mentioned
at the end of Section~\ref{sec:HMC}. More precisely, we compute over a
trajectory the rate of rejected moves in the Step~4 of
Algorithm~\ref{algo:HMC}. The logarithmic plot in Figure~\ref{fi:rejection}
shows a linear fit with a slope of about $3.1$, which confirms the expected
scaling in~$\Dt^3$.

\begin{figure}[p]
  \centering
  \includegraphics[width=.8\textwidth]{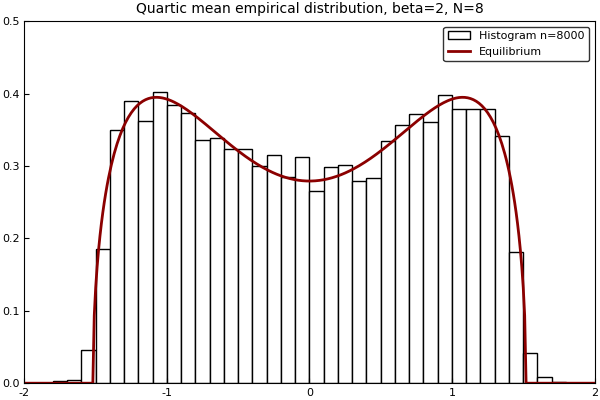}\\
  \includegraphics[width=.8\textwidth]{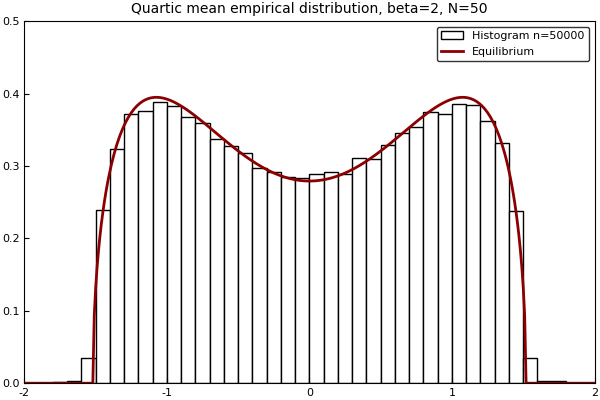}  
  \caption{\label{fi:quartic} Study of the quartic confinement with $N=8$
    (top) and $N=50$ (bottom). The solid line is the plot of the limiting
    spectral distribution~\eqref{eq:SC}. The bars form the histogram of
    simulations obtained using our HMC algorithm. This algorithm was run once
    with final-time $T=10^6$ and time-step $\Dt=0.5$. The histogram was
    produced by looking at the last half of the trajectory and retaining the
    positions each~$1000$ time-steps, producing $n$ values namely
    $\approx~8\times 10^3$ and $\approx~5\times 10^4$ respectively. We do not
    have a formula for the mean empirical distribution for this model. This
    gas describes the law of the eigenvalues of a random symmetric tridiagonal
    matrix model but its entries are not independent, see \cite[Prop.\
    2]{MR3433632}.}
\end{figure}

We also study the quartic confinement potential
$V(x)=x^4/4$, as in~\cite{MR3124976}. In this case, the empirical spectral
distribution is not known, but the equilibrium distribution has density with
respect to the Lebesgue measure given by
\[
  x\in\dR \mapsto(2a^2+x^2)\frac{\sqrt{4a^2-x^2}}{2\pi}\IND_{x\in[-2a,2a]},
  \quad a = 3^{-\frac{1}{4}}.
\]
The results of the numerical simulations, see Figure~\ref{fi:quartic}, show a
good agreement with the equilibrium measure when~$N$ is large. Note that a
tridiagonal random matrix model is known but it does not have independent
entries, see \cite[Prop.~2.1]{MR3433632}.

\subsection{Case study: 2D}
\label{sec:Ginibre}
\begin{figure}[ht!]
  \centering
  \includegraphics[width=.8\textwidth]{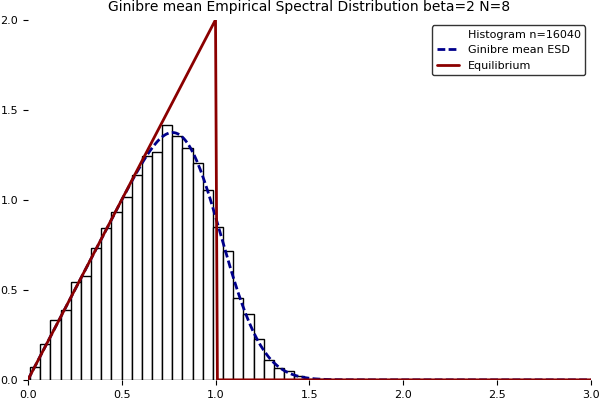}\\[2em]
  \includegraphics[width=.8\textwidth]{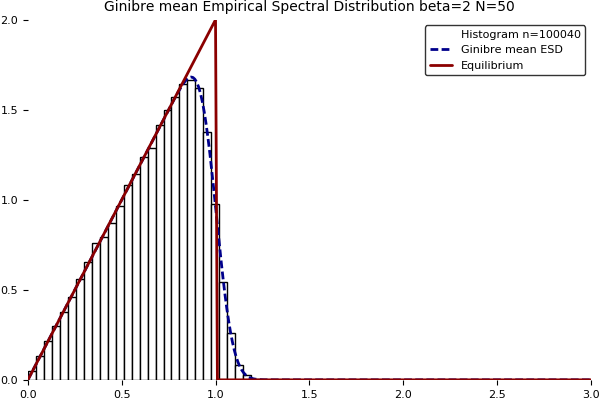}  
  \caption{\label{fi:Ginibre} Study of the complex Ginibre ensemble with $N=8$
    (top) and $N=50$ (bottom). The solid line is the plot of the limiting
    spectral distribution~\eqref{eq:G} while the dashed line is the plot of
    the mean empirical distribution~\eqref{eq:EGN}, both as functions of the
    radius $|z|$ and scaled by $2\pi$ (in order to obtain a radial density).
    The bars form the histogram of simulations obtained using our HMC
    algorithm. This algorithm was run $40$ times with final-time $T=10^5$ and
    time-step $\Dt=0.1$. The histogram was produced by looking at the last
    halves of the $40$ trajectories and retaining the positions each~$10000$
    time-steps, producing $n$ values namely $\approx~16\times 10^3$ and
    $\approx~10^5$ respectively.}
\end{figure}

\begin{figure}[p]
  \centering
  \includegraphics[width=.8\textwidth]{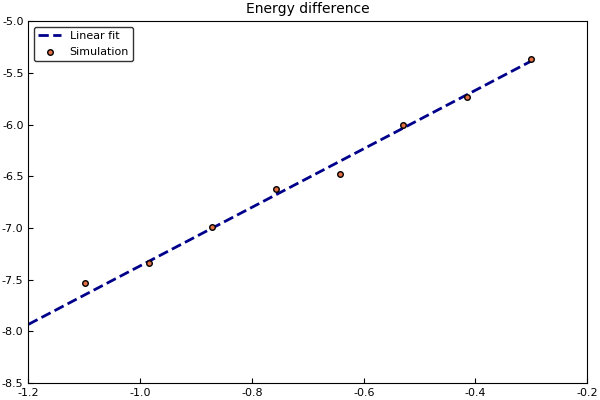}\\
  \caption{\label{fi:energy} Evolution of the energy difference in
    Algorithm~\ref{algo:HMC} as~$\Dt$ goes to zero, for the Complex Ginibre
    Ensemble with~$N=50$, $\beta=2$ and $T=10^{3}$ (in log-log coordinate). }
\end{figure}

We next consider in Figure~\ref{fi:Ginibre} the mean empirical distribution in
the case of the Complex Ginibre Ensemble \eqref{eq:BG} with $\be=2$, $N=8$. In
this case, we also know a theoretical formula for $\dE\mu_N$ under $P_N$,
given by~\eqref{eq:EGN}. For completeness, we investigate the scaling of the
relative energy difference in the Step~3 of Algorithm~\ref{algo:HMC} (by turning off
the selection procedure of Step~4). The logarithmic plot in Figure~\ref{fi:energy}
shows a slope of about $2.9$, which confirms the expected scaling in~$\Dt^3$ that
corresponds to the error of energy conservation, over one time step,
of the Verlet integrator~\eqref{eq:verlet}.

\begin{figure}[ht!]
  \centering
  \includegraphics[width=.6\textwidth]{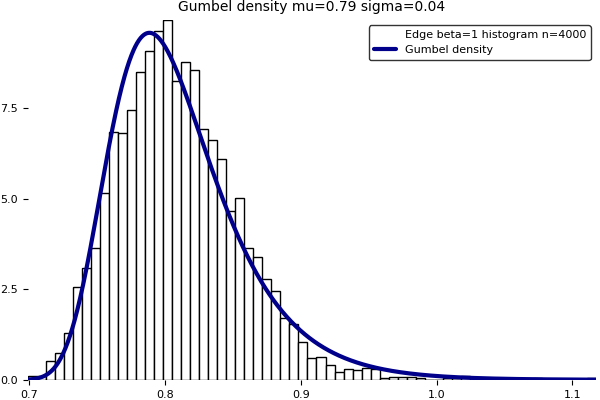}\\[2em]
  \includegraphics[width=.6\textwidth]{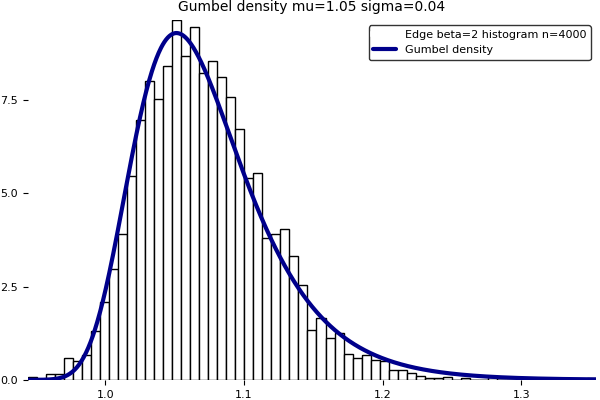}\\[2em]
  \includegraphics[width=.6\textwidth]{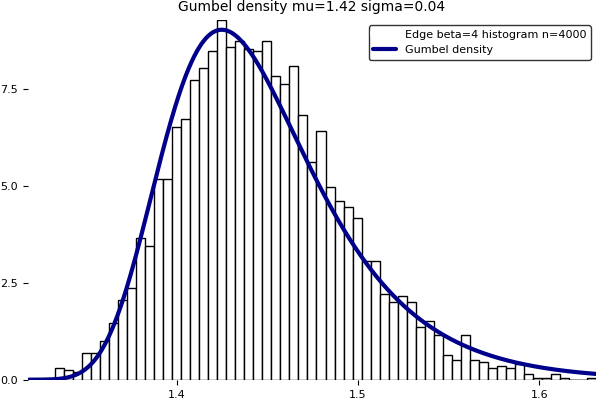} 
  \caption{\label{fi:GinibreEdge} Study of the fluctuation of the largest
    particle in modulus for the $\beta$ complex Ginibre ensemble with $N=50$,
    in the cases $\beta\in\{1,2,4\}$. The solid line is the plot of the fit
    with a translation-scale Gumbel distribution. The Gumbel fluctuation is
    proved only in the case $\beta=2$, see \cite{MR1986426,MR3215627}. These
    simulations suggest to conjecture that the Gumbel fluctuation is valid for
    any $\beta>0$. The simulation matches pretty well the edge support at
    $\sqrt{\beta/2}$ and suggests that the variance is not very sensitive to
    $\beta$.}
\end{figure}

We explore next in Figure~\ref{fi:GinibreEdge} the Gumbel fluctuation at the
edge, which is proved for $\be=2$ and conjectured for $\be\neq2$, see
\cite{MR1986426,MR3215627,dubach} (note that in this case we have a formula
for $\mu_*$ but not for $\mathbb{E}\mu_N$ under $P_N$). One could also explore
the crystallization phenomenon, see \cite{MR3429164} and references therein.

\subsection{Case study: 3D}

\begin{figure}[p]
  \centering
  \includegraphics[width=.8\textwidth]{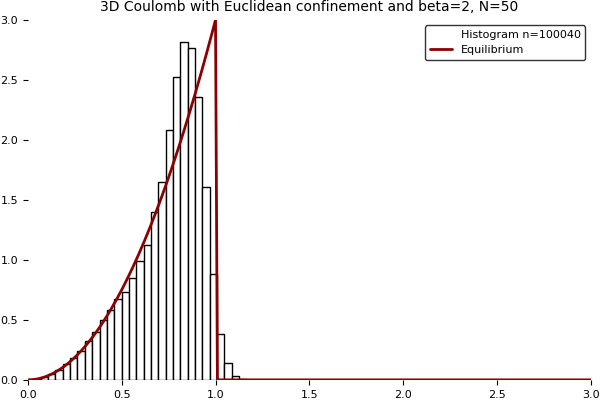}\\[2em]
  \includegraphics[width=.8\textwidth]{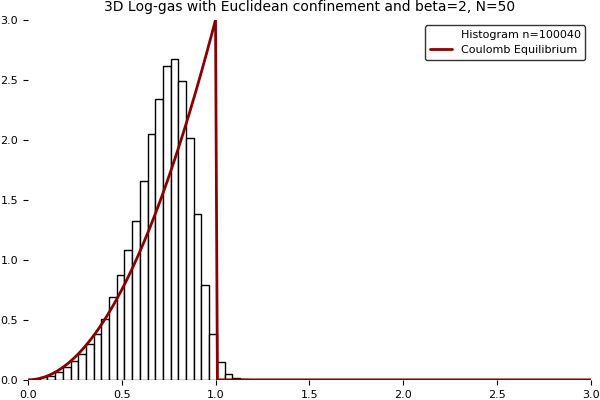}
  \caption{\label{fi:3d} Study of the 3D Coulomb case (top) and 3D Log-gas
    (bottom) with Euclidean confinement and $\beta=2$ and $N=50$. Equilibrium
    measure in solid line and histogram obtained with our HMC algorithm with
    $N=50$ and same simulation parameters as for Figure \ref{fi:Ginibre}. In
    contrast with the GUE case and the Ginibre case, we do not have a formula
    for the mean empirical distribution at fixed $N$ for both cases, and for
    the Log-gas (bottom) the equilibrium measure is not known.}
\end{figure}

In Figure \ref{fi:3d}, we finally turn to the Coulomb gas which corresponds to
$S=\dR^3$, $d=n=3$, $V=\ABS{\cdot}^2/\be$, $W=1/\ABS{\cdot}$ and to the
log-gas for which $W=-\log\ABS{\cdot}$. In the first case the equilibrium
measure $\mu_*$ is uniform on the centered ball of $\dR^d$ of radius
$(\be(d-2)/2)^{1/d}$, see for instance \cite[Cor.~1.3]{MR3262506}, while in
the second case the equilibrium measure is not know yet, see however
\cite{chafai-saff}. In both cases we do not have a formula for
$\mathbb{E}\mu_N$ under $P_N$. One could study the fluctuation at the edge,
which is conjectured to be Gumbel, just like for the complex Ginibre ensemble
in 2D.

\appendix

\section{Julia code}

Here is a program written in the Julia 
language\footnote{\url{http://JuliaLang.org/}} illustrating our method. 
It allows to exploit the multiple cores of modern processors and works in
parallel on clusters. Beware that this code is not fully optimized, for
instance the energy and its gradient could be computed simultaneously for
better performance.

\medskip

\begin{lstlisting}[basicstyle=\tiny]
#------------------------------------------------------------------------------#
#------------- Simulating coulomb gases with HMC algorithm --------------------#
#------------------------------------------------------------------------------#

# Tested with Julia 1.0. D. Chafai + G. Ferre : https://arxiv.org/abs/1806.05985

using Distributed # for @everywhere and nprocs()
@everywhere using Printf # for @sprintf()
@everywhere using LinearAlgebra # for norm()
@everywhere using DelimitedFiles # for Base.writedlm()

@everywhere begin # for parallel computing: julia -p NumberOfAdditionalProcesses
    #--------------------------------------------------------------------------#
    # Customization part : parameters, confinement, and interaction            #
    # -------------------------------------------------------------------------#

    ## Parameters. Note that in this code U_N(y)=|y|^2/2.
    
    # Final time and time step
    const T = 1e4
    const dt = 0.1
    # Number of eigenvalues
    const N = 8
    # Dimension of the physical space
    const dim = 1 # works for dimensions 1, 2, 3
    # Temperature and friction
    const beta = 2.
    # Riesz parameter for Riesz interaction
    const s = 1.

    ## Functions
    
    # Confinement potential V and its gradient
    @inline function confinement(x)
        return dot(x,x)/(2*beta) # 1D Beta-Hermite
        # return dot(x,x)/beta  # 2D Beta-Ginibre, 3D Beta-Coulomb  
    end
    @inline function confinement_gradient(x)
        return x/beta     # 1D Beta-Hermite
        # return 2*x/beta # 2D Beta-Ginibre, 3D Beta-Coulomb
    end
    
    # Interaction potential W and its gradient
    @inline function interaction(x,y)
        return -log(norm(x-y)) # 1D Beta-H., 2D Beta-Gin., 2D/3D Beta log-gas.
        # return 1/norm(x-y)   # 3D Beta-Coulomb 
        # return 1/norm(x-y)^s # Riesz
    end    
    @inline function interaction_gradient(x,y) 
        v = x-y
        return -v/norm(v)^2 # 1D Beta-H., 2D Beta-Gin        
        # return -v/norm(v)^3  # 3D Beta-Coulomb 
        # return -s*x/norm(v)^(s+2) # Riesz
    end

    #------------------------------------------------------------------#
    #--- Parameters computed from inputs ------------------------------#
    #------------------------------------------------------------------#

    const alphan = 1.
    const betan = beta * N^2
    const gamman = 1. / alphan
    # Parameters for discretisation of fluctuation-dissipation part L2
    const etan = exp(- gamman * alphan * dt)
    const sdn = sqrt((1-etan^2)/betan)
    #-- I/O parameter, write the configuration every niterio steps
    const niterio = 1000
    # Number of iterations and number of outputs
    const niter = Int64(round(T/dt))
    const nsteps = Int64(round(niter/niterio))
   
    #------------------------------------------------------------------#
    #---------- Core part - Be careful and good luck! -----------------#
    #------------------------------------------------------------------#
    
    ## Functions
    
    # Potential energy H_N
    @inline function energy(X)
        ener = 0
        @inbounds for i = 1:N
            @inbounds for j = i+1:N
                ener += interaction(X[i],X[j]) /N 
            end            
            ener += confinement(X[i]) 
        end
        return ener /N
    end # function energy()

    # Kinetic energy U_N
    @inline function kinetic(Y) 
        return norm(Y)^2 /2.
    end # function kinetic()

    # Force applied on particle at X[i] from all others at positions X[j] j!=i
    @inline function compute_force!(X,F) # -Grad H_N 
        # Computation of interaction forces between each pairs
        Fpairs = Array{Vector{Float64}}(undef, N,N) # we use only N(N-1)/2 entries
        @inbounds for i = 1:N
            @inbounds for j = 1:i-1
                Fpairs[i,j] = -interaction_gradient(X[i],X[j])
            end               
        end
        # Computation of total force on each particle       
        @inbounds for i = 1:N
            F[i] = zeros(dim)     
            # Interaction
            @inbounds for j = 1:i-1
                F[i] += Fpairs[i,j] 
            end               
            @inbounds for j = i+1:N
                F[i] -= Fpairs[j,i]
            end
            F[i] /= N
            # Confinement        
            F[i] -= confinement_gradient(X[i])
            F[i] /= N
        end
    end # function compute_force!()

    # compute the new force and speed
    @inline function verlet_integrator!(Fnew, Fcur, Xnew, Ynew, X, Y)
        @inbounds for i=1:N
            Ynew[i] = Y[i] + Fcur[i] * alphan * dt/2. 
            Xnew[i] = X[i] + Ynew[i] * alphan * dt 
        end
        compute_force!(Xnew,Fnew)
        @inbounds for i=1:N
            Ynew[i] += Fnew[i] * alphan * dt/2
        end
    end # function verlet_integrator!()

    # update positions and speed
    function update!(X, Y, Fcur, Xnew, Ynew, Fnew, Epot, acceptrate)
        #--- Speed resampling
        @inbounds for i = 1:N
            Y[i] = etan * Y[i] + sdn * randn(dim)
        end
        Ekin = kinetic(Y)
        Energy = Epot + Ekin
        #--- Verlet integrator. Position-speed proposal will be in (Xnew,Ynew).        
        verlet_integrator!(Fnew, Fcur, Xnew, Ynew, X, Y)
        # New energy
        Epotnew = energy(Xnew)
        Ekinnew = kinetic(Ynew)
        NewEnergy = Epotnew + Ekinnew
        # Metropolis ratio
        r = betan * (- NewEnergy + Energy)
        # Selection-rejection step
        if log(rand()) <= r
            # acceptation
            @inbounds @simd for i = 1:N
                X[i] = Xnew[i]
                Y[i] = Ynew[i]
                Fcur[i] = Fnew[i]
            end
            acceptrate[1] += 1
            Epot = Epotnew
        else # rejection: speed inversion    
            @inbounds @simd for i = 1:N
                Y[i] = -Y[i]
            end       
        end
        return Epot
    end # function update()

    # Runs a trajectory of HMC algorithm and compute averages
    function HMC(runid)
        #--- For output : for positions/velocities every niterio steps
        TrajectoryX = Array{Float64}(undef, nsteps, N*dim)
        TrajectoryY = Array{Float64}(undef, nsteps, N*dim)       
        #--- For output : Acceptation rate for the HMC selection step
        acceptrate = zeros(1)
        # Local variables
        #--- configuration and speed
        X = Vector{Vector{Float64}}(undef, N)
        Y = Vector{Vector{Float64}}(undef, N)
        #--- initial forces
        Fcur = Vector{Vector{Float64}}(undef, N)
        #--- Same quantities for the proposal
        Ynew = Vector{Vector{Float64}}(undef, N)
        Xnew = Vector{Vector{Float64}}(undef, N)
        Fnew = Vector{Vector{Float64}}(undef, N)
        # random initial configuration with uniform law on a square
        for i = 1:N
            X[i] = -1 .+ 2 * rand(dim)
        end
        if dim == 1
            X = sort(X)
        end
        # initial zero speed and forces
        for i = 1:N
            Y[i] = zeros(dim)
            Fcur[i] = zeros(dim)
            Xnew[i] = X[i]
            Ynew[i] = Y[i]
            Fnew[i] = Fcur[i]
        end
        #--- initialization of different quantities
        Ekin = kinetic(Y)
        Epot = energy(X)
        Energy = Epot + Ekin
        #--- Loop over time        
        @fastmath @inbounds for n = 1:niter            
            #--- save configuration every niterio steps
            if n % niterio == 0
                for i = 1:N
                    for k = 1:dim
                        TrajectoryX[ Int64(n/niterio), (i-1) * dim + k ] = X[i][k]
                        TrajectoryY[ Int64(n/niterio), (i-1) * dim + k ] = Y[i][k]                        
                    end
                end
            end
            #--- update positions and speeds
            Epot = update!( X, Y, Fcur, Xnew, Ynew, Fnew, Epot, acceptrate)
            
        end
        ## Post-processing
        print("Percentage of rejected steps: ", 1 .- acceptrate/niter, "\n")
        # Write the data in text files - Whole trajectory sample
        writedlm(@sprintf("positions-%i",runid), TrajectoryX, " ")
        writedlm(@sprintf("velocities-%i",runid), TrajectoryY, " ") 
    end # function HMC()
#
end # @everywhere

### Main part - runs only on main Julia process.
Nprocs = nprocs()
print("Number of processes is ",Nprocs,".\n")
print("Time of the simulation is ", T, ".\n")
print("Number of time steps is ", T/dt ,".\n")
## Launching computations on Nprocs parallel processes.
output = @time pmap(HMC,1:Nprocs)
### EOF
    
\end{lstlisting}

\subsection*{Acknowledgments} %
We warmly thank Gabriel Stoltz for his encouragements and for very useful
discussions on the theoretical and numerical sides of this work. We are also
grateful to Thomas Leblé and Laure Dumaz for their comments on the first
version.

\renewcommand{\MR}[1]{}


\begin{thebibliography}{10}

\bibitem{MR2760897}
{\scshape G.~W. Anderson, A.~Guionnet {\normalfont \smfandname} O.~Zeitouni} --
  \emph{An {I}ntroduction to {R}andom {M}atrices}, Cambridge Studies in
  Advanced Mathematics, vol. 118, Cambridge University Press, Cambridge, 2010.

\bibitem{bardenet-hardy}
{\scshape R.~Bardenet {\normalfont \smfandname} A.~Hardy} -- {\og {M}onte
  {C}arlo with {D}eterminantal {P}oint {P}rocesses\fg}, preprint
  \href{http://arxiv.org/abs/1605.00361v1}{arXiv:1605.00361v1}, 2016.

\bibitem{berman}
{\scshape R.~J. Berman} -- {\og On large deviations for {G}ibbs measures, mean
  energy and {G}amma-convergence\fg}, preprint
  \href{https://arxiv.org/abs/1610.08219v1}{arXiv:1610.08219v1}, 2016.

\bibitem{MR3129023}
{\scshape A.~Beskos, N.~Pillai, G.~Roberts, J.-M. Sanz-Serna {\normalfont
  \smfandname} A.~Stuart} -- {\og Optimal tuning of the hybrid {M}onte {C}arlo
  algorithm\fg}, \emph{Bernoulli} \textbf{19} (2013), no.~5A, p.~1501--1534.

\bibitem{MR3429164}
{\scshape X.~Blanc {\normalfont \smfandname} M.~Lewin} -- {\og The
  crystallization conjecture: a review\fg}, \emph{EMS Surv. Math. Sci.}
  \textbf{2} (2015), no.~2, p.~225--306.

\bibitem{coulsim}
{\scshape F.~Bolley, D.~Chafaï {\normalfont \smfandname} J.~Fontbona} -- {\og
  Dynamics of a planar {C}oulomb gas\fg}, preprint
  \href{http://arxiv.org/abs/1706.08776v3}{arXiv:1706.08776v3} to appear in
  {A}nnals of {A}pplied {P}robability, 2017.

\bibitem{bourabee-eberle-zimmer}
{\scshape N.~Bou-Rabee, A.~Eberle {\normalfont \smfandname} R.~Zimmer} -- {\og
  Coupling and {C}onvergence for {H}amiltonian {M}onte {C}arlo\fg}, preprint
  \href{http://arxiv.org/abs/1805.00452}{arXiv:1805.00452v1}, 2018.

\bibitem{bou2017geometric}
{\scshape N.~Bou-Rabee {\normalfont \smfandname} J.~M. Sanz-Serna} -- {\og
  Geometric integrators and the {H}amiltonian {M}onte {C}arlo method\fg},
  preprint \href{http://arxiv.org/abs/1711.05337}{arXiv:1711.05337v1}, 2017.

\bibitem{bouchard2018bouncy}
{\scshape A.~Bouchard-C{\^o}t{\'e}, S.~J. Vollmer {\normalfont \smfandname}
  A.~Doucet} -- {\og The bouncy particle sampler: A nonreversible
  rejection-free {M}arkov chain {M}onte {C}arlo method\fg}, \emph{J. Am. Stat.
  Assoc.} \textbf{113} (2018), no.~522, p.~855--867.

\bibitem{MR2742422}
{\scshape S.~Brooks, A.~Gelman, G.~L. Jones {\normalfont \smfandname} X.-L.
  Meng} (\smfedsname) -- \emph{Handbook of {M}arkov {C}hain {M}onte {C}arlo},
  Chapman \& Hall/CRC Handbooks of Modern Statistical Methods, CRC Press, Boca
  Raton, FL, 2011.

\bibitem{brosse2017tamed}
{\scshape N.~Brosse, A.~Durmus, {\'E}.~Moulines {\normalfont \smfandname}
  S.~Sabanis} -- {\og The tamed unadjusted {L}angevin algorithm\fg}, preprint
  \href{http://arXiv.org/abs/1710.05559v2}{arXiv:1710.05559v2}, 2017.

\bibitem{MR3262506}
{\scshape D.~Chafa{\"{\i}}, N.~Gozlan {\normalfont \smfandname} P.-A. Zitt} --
  {\og First-order global asymptotics for confined particles with singular pair
  repulsion\fg}, \emph{Ann. Appl. Probab.} \textbf{24} (2014), no.~6,
  p.~2371--2413.

\bibitem{MR3215627}
{\scshape D.~Chafa{\"{\i}} {\normalfont \smfandname} S.~P\'ech\'e} -- {\og A
  note on the second order universality at the edge of {C}oulomb gases on the
  plane\fg}, \emph{J. Stat. Phys.} \textbf{156} (2014), no.~2, p.~368--383.

\bibitem{chafai-saff}
{\scshape D.~Chafa{\"{\i}} {\normalfont \smfandname} E.~Saff} -- {\og Aspects
  of an {E}uclidean log-gas\fg}, work in progress, 2018.

\bibitem{coco}
{\scshape D.~Chafaï, A.~Hardy {\normalfont \smfandname} M.~Maïda} -- {\og
  Concentration for {C}oulomb gases and {C}oulomb transport inequalities\fg},
  preprint \href{https://arxiv.org/abs/1610.00980v3}{arXiv:1610.00980v3} to
  appear in J. Funct. Anal., 2018.

\bibitem{chafai-lehec}
{\scshape D.~Chafaï {\normalfont \smfandname} J.~Lehec} -- {\og On {P}oincaré
  and logarithmic {S}obolev inequalities for a class of singular {G}ibbs
  measures\fg}, preprint
  \href{http://arxiv.org/abs/1805.00708v2}{arXiv:1805.00708v2}, 2018.

\bibitem{dalalyan-riou-durand}
{\scshape A.~Dalalyan {\normalfont \smfandname} L.~Riou-Durand} -- {\og On
  sampling from a log-concave density using kinetic {L}angevin diffusions\fg},
  preprint \href{https://arxiv.org/abs/1807.09382}{arXiv:1807.09382v1}, 2018.

\bibitem{MR3439168}
{\scshape L.~Decreusefond, I.~Flint {\normalfont \smfandname} A.~Vergne} --
  {\og A note on the simulation of the {G}inibre point process\fg}, \emph{J.
  Appl. Probab.} \textbf{52} (2015), no.~4, p.~1003--1012.

\bibitem{duane-hmc}
{\scshape S.~Duane, A.~Kennedy, B.~J. Pendleton {\normalfont \smfandname}
  D.~Roweth} -- {\og {H}ybrid {M}onte {C}arlo\fg}, \emph{Physics Letters B}
  \textbf{195} (1987), no.~2, p.~216--222.

\bibitem{dubach}
{\scshape G.~Dubach} -- {\og Powers of {G}inibre {E}igenvalues\fg}, preprint
  \href{http://arxiv.org/abs/1711.03151}{arXiv:1711.03151v2}, 2017.

\bibitem{MR1936554}
{\scshape I.~Dumitriu {\normalfont \smfandname} A.~Edelman} -- {\og Matrix
  models for beta ensembles\fg}, \emph{J. Math. Phys.} \textbf{43} (2002),
  no.~11, p.~5830--5847.

\bibitem{2016JSP...163..457D}
{\scshape A.~B. {Duncan}, T.~{Leli{\`e}vre} {\normalfont \smfandname} G.~A.
  {Pavliotis}} -- {\og {Variance Reduction Using Nonreversible Langevin
  Samplers}\fg}, \emph{J. Stat. Phys.} \textbf{163} (2016), p.~457--491.

\bibitem{durmus2017convergence}
{\scshape A.~Durmus, E.~Moulines {\normalfont \smfandname} E.~Saksman} -- {\og
  On the convergence of {H}amiltonian {M}onte {C}arlo\fg}, preprint
  \href{http://arxiv.org/abs/1705.00166}{arXiv:1705.00166v1}, 2017.

\bibitem{MR2168344}
{\scshape A.~Edelman {\normalfont \smfandname} N.~R. Rao} -- {\og Random matrix
  theory\fg}, \emph{Acta Numer.} \textbf{14} (2005), p.~233--297.

\bibitem{MR3699468}
{\scshape L.~Erd{\H{o}}s {\normalfont \smfandname} H.-T. Yau} -- \emph{A
  dynamical approach to random matrix theory}, Courant Lecture Notes in
  Mathematics, vol.~28, Courant Institute of Mathematical Sciences, New York;
  American Mathematical Society, Providence, RI, 2017.

\bibitem{MR2377026}
{\scshape Z.~F. Ezawa} -- \emph{Quantum {H}all effects}, second \smfedname,
  World Scientific Publishing Co. Pte. Ltd., Hackensack, NJ, 2008, Field
  theoretical approach and related topics.

\bibitem{fathi2015error}
{\scshape M.~Fathi, A.-A. Homman {\normalfont \smfandname} G.~Stoltz} -- {\og
  Error analysis of the transport properties of {M}etropolized schemes\fg},
  \emph{ESAIM: Proceedings and Surveys} \textbf{48} (2015), p.~341--363.

\bibitem{MR2641363}
{\scshape P.~J. Forrester} -- \emph{Log-gases and {R}andom {M}atrices}, London
  Mathematical Society Monographs Series, vol.~34, Princeton University Press,
  Princeton, NJ, 2010.

\bibitem{MR3458536}
{\scshape P.~J. Forrester} -- {\og Analogies between random matrix ensembles
  and the one-component plasma in two-dimensions\fg}, \emph{Nuclear Phys. B}
  \textbf{904} (2016), p.~253--281.

\bibitem{garcia-zelada}
{\scshape D.~García-Zelada} -- {\og A large deviation principle for empirical
  measures on {P}olish spaces: {A}pplication to singular {G}ibbs measures on
  manifolds\fg}, preprint
  \href{https://arxiv.org/abs/1703.02680v2}{arXiv:1703.02680v2}, 2017.

\bibitem{hairer2006geometric}
{\scshape E.~Hairer, C.~Lubich {\normalfont \smfandname} G.~Wanner} --
  \emph{Geometric {N}umerical {I}ntegration: {S}tructure-{P}reserving
  {A}lgorithms for {O}rdinary {D}ifferential {E}quations}, Springer Series in
  Computational Mathematics, vol.~31, Springer Science \& Business Media, 2006.

\bibitem{hardy}
{\scshape A.~Hardy} -- {\og Polynomial {E}nsembles and {R}ecurrence
  {C}oefficients\fg}, preprint
  \href{http://arXiv.org/abs/1709.01287}{arXiv:1709.01287v1}, 2017.

\bibitem{MR3308615}
{\scshape L.~L. Helms} -- \emph{Potential {T}heory}, second \smfedname,
  Universitext, Springer, London, 2014.

\bibitem{MR3214779}
{\scshape M.~D. Hoffman {\normalfont \smfandname} A.~Gelman} -- {\og The
  {N}o-{U}-{T}urn sampler: adaptively setting path lengths in {H}amiltonian
  {M}onte {C}arlo\fg}, \emph{J. Mach. Learn. Res.} \textbf{15} (2014),
  p.~1593--1623.

\bibitem{MR3659634}
{\scshape T.~A. H\"oft {\normalfont \smfandname} B.~K. Alpert} -- {\og Fast
  updating multipole {C}oulombic potential calculation\fg}, \emph{SIAM J. Sci.
  Comput.} \textbf{39} (2017), no.~3, p.~A1038--A1061.

\bibitem{horowitz1991generalized}
{\scshape A.~M. Horowitz} -- {\og A generalized guided {M}onte {C}arlo
  algorithm\fg}, \emph{Phys. Lett. B} \textbf{268} (1991), no.~CERN-TH-6172-91,
  p.~247--252.

\bibitem{MR2216966}
{\scshape J.~B. Hough, M.~Krishnapur, Y.~Peres {\normalfont \smfandname}
  B.~Vir\'ag} -- {\og Determinantal processes and independence\fg},
  \emph{Probab. Surv.} \textbf{3} (2006), p.~206--229.

\bibitem{hutzenthaler2012strong}
{\scshape M.~Hutzenthaler, A.~Jentzen {\normalfont \smfandname} P.~E. Kloeden}
  -- {\og Strong convergence of an explicit numerical method for {SDE}s with
  nonglobally {L}ipschitz continuous coefficients\fg}, \emph{Ann. Appl.
  Probab.} \textbf{22} (2012), no.~4, p.~1611--1641.

\bibitem{MR3615091}
{\scshape T.~Jiang {\normalfont \smfandname} Y.~Qi} -- {\og Spectral radii of
  large non-{H}ermitian random matrices\fg}, \emph{J. Theoret. Probab.}
  \textbf{30} (2017), no.~1, p.~326--364.

\bibitem{jones2011adaptive}
{\scshape A.~Jones {\normalfont \smfandname} B.~Leimkuhler} -- {\og Adaptive
  stochastic methods for sampling driven molecular systems\fg}, \emph{J. Chem.
  Phys.} \textbf{135} (2011), no.~8, p.~084125.

\bibitem{Kapfer-2016}
{\scshape S.~C. Kapfer {\normalfont \smfandname} W.~Krauth} -- {\og Cell-veto
  {M}onte {C}arlo algorithm for long-range systems\fg}, \emph{Physical Review
  E} \textbf{94} (2016).

\bibitem{MR1260431}
{\scshape P.~E. Kloeden, E.~Platen {\normalfont \smfandname} H.~Schurz} --
  \emph{Numerical solution of {S}{D}{E} through computer experiments},
  Universitext, Springer-Verlag, Berlin, 1994, With 1 IBM-PC floppy disk (3.5
  inch; HD).

\bibitem{MR3433632}
{\scshape M.~Krishnapur, B.~Rider {\normalfont \smfandname} B.~Vir\'ag} -- {\og
  Universality of the stochastic {A}iry operator\fg}, \emph{Comm. Pure Appl.
  Math.} \textbf{69} (2016), no.~1, p.~145--199.

\bibitem{MR0350027}
{\scshape N.~S. Landkof} -- \emph{Foundations of {M}odern {P}otential
  {T}heory}, Springer-Verlag, New York-Heidelberg, 1972, Translated from the
  Russian by A. P. Doohovskoy, Die Grundlehren der mathematischen
  Wissenschaften, Band 180.

\bibitem{MR3382600}
{\scshape F.~Lavancier, J.~M{\o}ller {\normalfont \smfandname} E.~Rubak} --
  {\og Determinantal point process models and statistical inference\fg},
  \emph{J. R. Stat. Soc. Ser. B. Stat. Methodol.} \textbf{77} (2015), no.~4,
  p.~853--877.

\bibitem{MR2041832}
{\scshape M.~Ledoux} -- {\og Differential operators and spectral distributions
  of invariant ensembles from the classical orthogonal polynomials. {T}he
  continuous case\fg}, \emph{Electron. J. Probab.} \textbf{9} (2004), p.~no. 7,
  177--208.

\bibitem{lee-vempala}
{\scshape Y.~T. Lee {\normalfont \smfandname} S.~S. Vempala} -- {\og
  Convergence {R}ate of {R}iemannian {H}amiltonian {M}onte {C}arlo and {F}aster
  {P}olytope {V}olume {C}omputation\fg}, preprint
  \href{http://arxiv.org/abs/1710.06261}{arXiv:1710.06261v1}, 2017.

\bibitem{leimkuhler2015computation}
{\scshape B.~Leimkuhler, C.~Matthews {\normalfont \smfandname} G.~Stoltz} --
  {\og The computation of averages from equilibrium and nonequilibrium
  {L}angevin molecular dynamics\fg}, \emph{IMA J. Numer. Anal.} \textbf{36}
  (2015), no.~1, p.~13--79.

\bibitem{2013JSP...152..237L}
{\scshape T.~{Leli{\`e}vre}, F.~{Nier} {\normalfont \smfandname} G.~A.
  {Pavliotis}} -- {\og {Optimal Non-reversible Linear Drift for the Convergence
  to Equilibrium of a Diffusion}\fg}, \emph{J. Stat. Phys.} \textbf{152}
  (2013), p.~237--274.

\bibitem{MR2681239}
{\scshape T.~Leli\`evre, M.~Rousset {\normalfont \smfandname} G.~Stoltz} --
  \emph{Free {E}nergy {C}omputations. {A} {M}athematical {P}erspective},
  Imperial College Press, London, 2010.

\bibitem{MR2945148}
\bysame , {\og Langevin dynamics with constraints and computation of free
  energy differences\fg}, \emph{Math. Comp.} \textbf{81} (2012), no.~280,
  p.~2071--2125.

\bibitem{MR3509213}
{\scshape T.~Leli\`evre {\normalfont \smfandname} G.~Stoltz} -- {\og Partial
  differential equations and stochastic methods in molecular dynamics\fg},
  \emph{Acta Numer.} \textbf{25} (2016), p.~681--880.

\bibitem{levesque-verlet}
{\scshape D.~Levesque {\normalfont \smfandname} L.~Verlet} -- {\og On the
  theory of classical fluids {II}\fg}, \emph{Physica} \textbf{28} (1962),
  no.~11, p.~1124--1142.

\bibitem{levesque-verlet-bis}
\bysame , {\og Computer "{E}xperiments" on {C}lassical {F}luids. {III}.
  {T}ime-{D}ependent {S}elf-{C}orrelation {F}unctions\fg}, \emph{Phys. Rev. A}
  \textbf{2} (1970), p.~2514.

\bibitem{MR3124976}
{\scshape X.~H. Li {\normalfont \smfandname} G.~Menon} -- {\og Numerical
  solution of {D}yson {B}rownian motion and a sampling scheme for invariant
  matrix ensembles\fg}, \emph{J. Stat. Phys.} \textbf{153} (2013), no.~5,
  p.~801--812.

\bibitem{MATTINGLY2002185}
{\scshape J.~Mattingly, A.~Stuart {\normalfont \smfandname} D.~Higham} -- {\og
  Ergodicity for {SDE}s and approximations: locally {L}ipschitz vector fields
  and degenerate noise\fg}, \emph{Stoch. Proc. Appl.} \textbf{101} (2002),
  no.~2, p.~185 -- 232.

\bibitem{MR2129906}
{\scshape M.~L. Mehta} -- \emph{Random {M}atrices}, third \smfedname, Pure and
  Applied Mathematics (Amsterdam), vol. 142, Elsevier/Academic Press,
  Amsterdam, 2004.

\bibitem{milstein2013stochastic}
{\scshape G.~N. Milstein {\normalfont \smfandname} M.~V. Tretyakov} --
  \emph{Stochastic {N}umerics for {M}athematical {P}hysics}, Springer Science
  \& Business Media, 2013.

\bibitem{MR3334666}
{\scshape S.~Olver, R.~R. Nadakuditi {\normalfont \smfandname} T.~Trogdon} --
  {\og Sampling unitary ensembles\fg}, \emph{Random Matrices Theory Appl.}
  \textbf{4} (2015), no.~1, p.~1550002, 22.

\bibitem{MR1986426}
{\scshape B.~Rider} -- {\og A limit theorem at the edge of a non-{H}ermitian
  random matrix ensemble\fg}, \emph{J. Phys. A} \textbf{36} (2003), no.~12,
  p.~3401--3409, Random matrix theory.

\bibitem{MR2080278}
{\scshape C.~P. Robert {\normalfont \smfandname} G.~Casella} -- \emph{Monte
  {C}arlo {S}tatistical {M}ethods}, second \smfedname, Springer Texts in
  Statistics, Springer-Verlag, New York, 2004.

\bibitem{roberts2001optimal}
{\scshape G.~O. Roberts {\normalfont \smfandname} J.~S. Rosenthal} -- {\og
  Optimal scaling for various {M}etropolis-{H}astings algorithms\fg},
  \emph{Statistical science} \textbf{16} (2001), no.~4, p.~351--367.

\bibitem{roberts1996exponential}
{\scshape G.~O. Roberts, R.~L. Tweedie {\normalfont et~al.}} -- {\og
  Exponential convergence of {L}angevin distributions and their discrete
  approximations\fg}, \emph{Bernoulli} \textbf{2} (1996), no.~4, p.~341--363.

\bibitem{rossky-doll-friedman}
{\scshape P.~H. Rossky, J.~D. Doll {\normalfont \smfandname} H.~L. Friedman} --
  {\og Brownian dynamics as smart {M}onte {C}arlo simulation\fg}, \emph{The
  Journal of Chemical Physics} \textbf{69} (1978), p.~4628.

\bibitem{MR1485778}
{\scshape E.~B. Saff {\normalfont \smfandname} V.~Totik} -- \emph{Logarithmic
  {P}otentials with {E}xternal {F}ields}, Grundlehren der Mathematischen
  Wissenschaften [Fundamental Principles of Mathematical Sciences], vol. 316,
  Springer-Verlag, Berlin, 1997, Appendix B by Thomas Bloom.

\bibitem{MR2551206}
{\scshape A.~Scardicchio, C.~E. Zachary {\normalfont \smfandname} S.~Torquato}
  -- {\og Statistical properties of determinantal point processes in
  high-dimensional {E}uclidean spaces\fg}, \emph{Phys. Rev. E (3)} \textbf{79}
  (2009), no.~4, p.~041108, 19.

\bibitem{MR3309890}
{\scshape S.~Serfaty} -- \emph{Coulomb gases and {G}inzburg-{L}andau vortices},
  Zurich Lectures in Advanced Mathematics, Euro. Math. Soc. (EMS), Z\"urich,
  2015.

\bibitem{serfaty-icm2018}
\bysame , {\og Systems of points with {C}oulomb interactions\fg}, preprint
  \href{https://arxiv.org/abs/1712.04095v1}{arXiv:1712.04095v1}, 2017.

\bibitem{MR1631413}
{\scshape S.~Smale} -- {\og Mathematical problems for the next century\fg},
  \emph{Math. Intelligencer} \textbf{20} (1998), no.~2, p.~7--15.

\bibitem{MR1754783}
\bysame , {\og Mathematical problems for the next century\fg}, in
  \emph{Mathematics: frontiers and perspectives}, Amer. Math. Soc., Providence,
  RI, 2000, p.~271--294.

\bibitem{stoltz-trstanova}
{\scshape G.~Stoltz {\normalfont \smfandname} Z.~Trstanova} -- {\og Stable and
  accurate schemes for {L}angevin dynamics with general kinetic energies\fg},
  preprint \href{http://arxiv.org/abs/1609.02891v1}{arXiv:1609.02891v1}, 2016.

\bibitem{PDMPs}
{\scshape P.~Vanetti, A.~Bouchard-Côté, G.~Deligiannidis {\normalfont
  \smfandname} A.~Doucet} -- {\og Piecewise-{D}eterministic {M}arkov {C}hain
  {M}onte {C}arlo\fg}, preprint
  \href{https://arxiv.org/abs/1707.05296}{arXiv:1707.05296v1}, 2018.

\bibitem{verlet}
{\scshape L.~Verlet} -- {\og Computer "{E}xperiments" on {C}lassical {F}luids.
  {I}. {T}hermodynamical {P}roperties of {L}ennard-{J}ones {M}olecules\fg},
  \emph{Phys. Rev.} \textbf{159} (1967), no.~98.

\end{thebibliography}
%

\providecommand{\bysame}{\leavevmode ---\ }
\providecommand{\og}{``}
\providecommand{\fg}{''}
\providecommand{\smfandname}{\&}
\providecommand{\smfedsname}{eds.}
\providecommand{\smfedname}{ed.}
\providecommand{\smfmastersthesisname}{Memoir}
\providecommand{\smfphdthesisname}{Thesis}

\end{document}